\documentclass{amsart}
\usepackage{amssymb,amsmath}
\usepackage{ytableau}
\usepackage[enableskew]{youngtab}
\usepackage{young}
\usepackage{tikz,pgf}

\newtheorem{theorem}{Theorem}

\newtheorem{conjecture}{Conjecture}

\begin{document}

\title{Restricted $k$-color partitions, II}
\author{William J. Keith}
\maketitle

\begin{abstract}

We consider $(k,j)$-colored partitions, partitions in which $k$ colors exist but at most $j$ colors may be chosen per size of part.  In particular these generalize overpartitions.  Advancing previous work, we find new congruences, including in previously unexplored cases where $k$ and $j$ are not coprime, as well as some noncongruences. As a useful aside, we give the apparently new generating function for the number of partitions in the $N \times M$ box with a given number of part sizes, and extend to multiple colors a conjecture of Dousse and Kim on unimodality in overpartitions.

\end{abstract}

\section{Introduction}

Two classes of partition-theoretic objects that have excited much research are overpartitions, defined in 2004 in \cite{CoLo} and since studied in \cite{DoKim}, \cite{HirschSell}, and \cite{XiaYao} among many other papers; and $k$-colored partitions or multipartitions, which have a longer history, being of interest to see Ramanujan \cite{BGK} and many authors since: see \cite{GEA3} for a good survey.  Precise definitions for these are given below.

In this paper, second in a series beginning in \cite{Keith1}, we further study a unification and generalization of these two classes, $(k,j)$-colored partitions.  Think of overpartitions as 2-colored partitions in which only one color may be used for each part size: the $(k,j)$-colored partitions are the natural generalization in which, of $k$ colors, only $j$ may be used for a given part size.  The usual $k$-colored partitions are $(k,k)$-colored partitions, and overpartitions are $(2,1)$-colored partitions.

Overpartitions exhibit many congruences and symmetries.  Many instances of the $(k,j)$-colored partitions with particular values of $k$ and $j$ were shown in \cite{Keith1} to possess such features.  Most of the pairs $(k,j)$ studied in that paper were coprime, since those were immediately accessible to more elementary machinery. Emphasis was given to $j=1$ and $j=k-1$.  Here we consider some non-coprime cases, which, though more intricate, exhibit new internal symmetries and connections to other functions such as divisor sums.  As with the previous paper, the congruences proved herein are certainly samples of a wide and rich field of similar results.

In addition to the above results we consider $(k,1)$-colored partitions in the $M \times N$ box, and generalize a conjecture of Dousse and Kim on unimodality of overpartitions.  In order to do so we require a theorem, apparently new, giving the generating function for the number of partitions in the $M \times N$ box in which exactly $r$ different part sizes appear.  This completes a trio of theorems on partitions with specified numbers of part sizes, with the unrestricted case by MacMahon \cite{MacMahon} and elaborated by Andrews \cite{GEA1}, and the singly-restricted case in which sizes must be less than $N$ given just in 2016 by Merca in \cite{Merca}.

The main theorems of the paper are as follows.

Denote by $c_{k,j}(n)$ the number of $(k,j)$-colored partitions of $n$.

\begin{theorem}\label{mod9} For $j \in \{2,5,8,9\}$, we have $c_{9,j}(3n+2) \equiv 0 \pmod{27}$ for all $n \geq 0$. For $j \in \{3,6\}$, we have $c_{9,j}(9n+2) \equiv 0 \pmod{27}$ for all $n \geq 0$. For $j \in \{1,4,7\}$, we have that if $c_{9,j}(3^\ell n+2) \not\equiv 0 \pmod{27}$ for all $n \geq 0$, then $\ell$ is large, and does not exist for $j=1$.
\end{theorem}

Denote by $\overline{p}(n)$ the overpartition function.  Then

\begin{theorem}\label{2pp} $$c_{2p,p}(n) \equiv \overline{p} \left( \frac{n}{p} \right) + p \sum_{\lambda \in S} 2^{r(\lambda)} (t(\lambda)+1) y(\lambda)^{-1} \pmod{p^2}$$ where $S$ is a certain subset of partitions of $n$, $t(\lambda)$, $y(\lambda)$ are statistics on these partitions, and inversion is taken mod $p$.
\end{theorem}

Given $\vec{t}=(t_1,\dots,t_k)$, let $p_{M \times N}(\vec{t};n)$ denote the number of $(k,1)$-colored partitions of $n$, i.e. $k$-colored partitions where only one color is chosen per size of part, with at most $N$ parts, all of size at most $M$, in which $t_i$ part sizes are of color $i$ for all $i$.  Then

\begin{theorem}\label{K1Box} Fix $M$ and $N$ dimensions of a box and fix a number of colors $k$.  We have the generating function $$\sum_{n,\vec{t}} p_{M \times N} (\vec{t};n) x_1^{t_1} \dots x_k^{t_k} q^n = \sum_{r=0}^{min(M,N)} (x_1 + \dots + x_k)^r f_{(M \times N; r)}(q),$$ where $\vec{t}$ runs over nonnegative integer vectors in ${\mathbb{Z}_{\geq 0}}^k$ and $f_{(M \times N; r)}(q)$ is the generating function for partitions in the $M \times N$ box with exactly $r$ different part sizes.
\end{theorem}

We record this latter generating function as being of significant independent interest:

\begin{theorem}\label{SizesInBox} $$f_{(M \times N);r}(q) = \sum_{r_0 = r}^{min(M,N)} \binom{r_0}{r} (-1)^{r_0-r} q^{\binom{r_0+1}{2}} \left[ {M \atop {r_0}} \right]_q \left[ {{M+N-r_0} \atop {N-r_0}} \right]_q.$$
\end{theorem}

Finally, we propose a generalization of a unimodality conjecture of Dousse and Kim on overpartitions in the box:

\begin{conjecture}\label{DKGen} Fix $n$.  Set any one $x_i = 1$ and consider the remaining $(n-1)$-dimensional array of coefficients in the $x_j$ on $q^n$ in $\sum_{n,\vec{t}} p_{M \times N} (t_1,\dots,t_k;n) x_1^{t_1} \dots x_k^{t_k} q^n$.  Any linear cut through this array is a unimodal sequence.
\end{conjecture}

This generalizes the conjecture of Dousse and Kim on overpartitions, which is the $k=2$ case in which $x_2$ is set to 1 and the remaining 1-dimensional sequence of coefficients is conjectured to be unimodal.

\section{Background and definitions}

The most common notations we will use for partitions are either the sum notation, $\lambda = \lambda_1 + \lambda _ 2 + \dots \vdash N$, with the $\lambda_i$ the parts of the partition; or the frequency notation, $\lambda = \lambda_1^{e_1} \lambda_2^{e_2} \dots $, meaning that part size $\lambda_1$ appears $e_1$ times, etc.  Context should always make it clear which notation is being used, as we never talk about exponents for parts.

Overpartitions are partitions in which the last instance of a given part may be overlined, or not.  The overpartitions of 3 are \[ 3, \overline{3}, 2+1, \overline{2}+1, 2+\overline{1}, \overline{2}+ \overline{1}, 1+1+1, 1+1+\overline{1} .\]

The generating function for $\overline{p}(n)$, the number of overpartitions of $n$, is \[ \overline{P}(q) := \sum_{n=0}^\infty \overline{p}(n) q^n = \prod_{k=1}^\infty \frac{1+q^k}{1-q^k} = \prod_{k=1}^\infty \frac{1-q^{2k}}{(1-q^k)^2}. \]

The $k$-colored partitions, also known as multipartitions, are partitions in which parts may appear in $k$ different types known as colors, with the order of colors not mattering (i.e., by convention we list colors within a part in weakly descending order).  For instance, denoting colors by subscripts, the 2-colored partitions of 3 are \[ 3_2, 3_1, 2_2+1_2, 2_2 + 1_1, 2_1+1_2, 2_1+1_1, 1_2+1_2+1_2, 1_2+1_2+1_1,1_2+1_1+1_1, 1_1+1_1+1_1 .\]

The generating function of $c_k(n)$, the number of $k$-colored partitions of $n$, is \[ C_k(q) := \sum_{n=0}^\infty c_k(n) q^n = \prod_{n=1}^\infty \frac{1}{(1-q^n)^k}. \]

In the language of colored partitions, overpartitions are 2-colored partitions with the restriction that only one color may be used per size of part.  This leads to the natural generalization in which, of $k$ colors, only $j$ may be used per size of part, which we call $(k,j)$-colored partitions.  Call $c_{k,j}(n)$ the number of $(k,j)$-colored partitions of $n$.  Their generating function is

\begin{multline*}C_{k,j}(q) := \sum_{n=0}^\infty c_{k,j}(n) q^n = \prod_{n=1}^\infty \left( 1 + \frac{\binom{k}{1} q^n}{1-q^n} + \frac{\binom{k}{2} q^{2n}}{(1-q^{n})^2} + \dots + \frac{\binom{k}{j}q^{jn}}{(1-q^n)^j} \right) \\
= \frac{1}{{(q)_\infty}^j} \prod_{n=1}^\infty \left( \sum_{i=0}^j \binom{k}{i} (1-q^n)^{j-i} q^{in} \right) = \frac{1}{{(q)_\infty}^j} \prod_{n=1}^\infty \left( \sum_{i=0}^j \binom{k-j+i-1}{i} q^{i n} \right).
\end{multline*}

Here we have used the standard notation \[ (a;q)_n = (1-a)(1-aq)\dots(1-aq^{n-1}) , \quad (a;q)_\infty = \lim_{n\rightarrow \infty} (a;q)_n, \quad (q)_\infty := (q;q)_\infty .\]

Partitions in which there are at most $N$ parts, each of size at most $M$, are said to fit in the $M \times N$ box, and have the $q$-binomial coefficients as their generating function:

\[ \left[ {{M+N} \atop N} \right]_q = \frac{(q)_{M+N}}{(q)_M(q)_N} . \]

P. A. MacMahon, and later George Andrews, studied the generating function for otherwise unrestricted partitions with exactly $r$ part sizes:

$$N_r(q) = \frac{1}{(q;q)_\infty} \sum_{n=0}^\infty (-1)^{n-r} \frac{\binom{n}{r}q^{\binom{n+1}{2}}}{(q;q)_n}.$$

Just in 2016, Mircea Merca in produced the generating function for partitions with $r$ part sizes restricted by largest part size:

$$N_{r,\leq M} (q) = \frac{1}{(q;q)_M} \sum_{j=r}^M (-1)^{j-r} q^{\binom{j+1}{2}} \binom{j}{r} \left[ {M \atop j} \right]_q.$$

Our Theorem \ref{SizesInBox} doubly restricts the partitions and is thus the completing third member in this line.  It is easily seen that Theorem \ref{SizesInBox} becomes Merca's identity in the limit $N \rightarrow \infty$.

Finally, when for two power series $f(q) = \sum_{n=0}^\infty a(n) q^n$ and $g(q) = \sum_{n=0}^\infty b(n) q^n$ we write $f(q) \equiv_m g(q)$, we mean $a(n) \equiv b(n) \pmod{m}$ for all $n$.

\section{Proofs of the main results}

\subsection{Proof of Theorem \ref{mod9}}  Begin by noticing that, given a partition $\lambda$, there are a large number of $(k,j)$-colored partitions which assign different colors to part sizes in $\lambda$.  Many of these can be grouped and ignored when proving congruences such as this theorem.  For instance, among all assignments of colors to a partition in which the part size 100 appears 5 times, those in which this part size is assigned two colors, the first coloring three parts and the next coloring 2, can be grouped in the $\binom{9}{2}$ ways to choose these colors.

Recall that for $i \not\equiv 0 \pmod{3}$, $\binom{9}{i} \equiv 0 \pmod{9}$, while $\binom{9}{3} \equiv \binom{9}{6} \equiv 3 \pmod{9}$.  Hence for $j < 9$, $(9,j)$-colored partitions with three or more distinct sizes of part contribute groups of size divisible by 27 to the total, and so we need only consider partitions in which 1 or 2 part sizes appear.

For $j=2$, we are considering 9-colored partitions in which only 1 or 2 colors may be used per size of part.

If a part size has 1 color, there are 9 options; if a part size has 2 colors, there are $\binom{9}{2}=36 \equiv 0 \pmod{9}$ options that assign the same number of parts to the first and second.  Ergo, if a partition $\lambda$ has two or more part sizes, then the associated $(9,2)$-colored partitions number 0 mod 81.

This reduces the problem to observing $\lambda = \lambda_1^{e_1} \vdash N$, partitions into just 1 part size, which contribute $$\sum_{e_1 \vert N} \binom{9}{1} \binom{e_1-1}{0} + \binom{9}{2}\binom{e_1-1}{1} \equiv 9 \sigma_1(N) \pmod{27}.$$

But for $N \equiv 2 \pmod{3}$, $N$ has at least one prime factor $p_i \equiv 2 \pmod{3}$ appearing to an odd power, so $\sigma_1(N) \equiv 0 \pmod{3}$, and so $9\sigma_1(N) \equiv 0 \pmod{27}$.

For $j=5$, partitions with two distinct part sizes contribute a nonzero amount mod 27 only if both part sizes include exactly 3 of the 9 possible colors.  Let $\lambda = \lambda_1^{e_1} \lambda_2^{e_2}$.  The number of $(9,5)$-colored partitions of this type contributed is $$\binom{9}{3} \binom{9}{3} \binom{e_1-1}{2} \binom{e_2-1}{2}.$$  But $\binom{k}{2} \equiv 0 \pmod{3}$ unless $k \equiv 2 \pmod{3}$.  If $e_1 - 1 \equiv e_2 - 1 \equiv 2 \pmod{3}$, then $e_1 \equiv e_2 \equiv 0 \pmod{3}$, which is impossible since $N \not\equiv 0 \pmod{3}$.  Hence partitions of this type contribute a number divisible by 27 to the total.

It remains to consider partitions into exactly one part size, which contribute a total of $$\sum_{e_1 \vert N} \sum_{i=1}^5 \binom{9}{i} \binom{e_1 - 1}{i-1} = \frac{3}{4} \left(76-162e_1+133{e_1}^2-42{e_1}^3+7{e_1}^4\right).$$

All divisors $e_1$ of $N \equiv 2 \pmod{3}$ are either 1 or 2 mod 3, and furthermore these come in pairs.  That is, if $e_1 \equiv 1 \pmod{3}$, the complementary divisor $N/e_1 \equiv 2 \pmod{3}$.  It is quickly seen that if $e_1 \equiv 1 \pmod{3}$, then the sum above is congruent to 9 mod 27, while if $e_1 \equiv 2 \pmod{3}$, the sum is 18 mod 27.  Hence the contributions of complementary divisors sum to 0 mod 27.

These exhaust all cases for $j=5$ and so the theorem is proved for this case.

The $j=8$ case is established with a similar argument; $\binom{9}{6}=\binom{9}{3}$, and the binomials $\binom{e_i-1}{5}$ behave the same way as $\binom{e_i-1}{2}$ mod 3.  The polynomial $\sum_{i=1}^8 \binom{9}{i}\binom{e_1-1}{i-1}$ has the same property for complementary divisors as discussed earlier.

Interestingly, the case $j=9$ can be proved from the case $j=8$: it was shown in \cite{Keith1} that the generating functions $C_{k,k-1}$ are $\eta$-quotients, and in fact we have \[ C_{9,8}(q) = \frac{(q^9;q^9)_\infty}{(q;q)_\infty^9} \quad , \quad C_{9,9}(q) = \frac{1}{(q;q)_\infty^9} . \]

This means that the coefficients $c_{9,9}(n)$ can be written with the (9-magnified) pentagonal number recurrence in terms of $c_{9,8}(n)$:

$$c_{9,9}(n) = c_{9,8}(n) - c_{9,8}(n-9) - c_{9,8}(n-18) + c_{9,8}(n-45) + c_{9,8}(n-63) - \dots .$$

If $N \equiv 2 \pmod{3}$, all terms $N-9t$ of the recurrence for $c_{9,9}(N)$ lie in the same arithmetic progression, and hence all are 0 mod 27.  This is perhaps a relatively rare case of a combinatorial proof of a congruence for a multipartition function above the number of colors itself.

For $j=3$ and $j=6$, we refer to the generating function.  $$C_{9,3}(q) \equiv_{27} \prod_{n=1}^\infty \frac{1+6q^n-6q^{2n}+2q^{3n}}{(1-q^n)^3} = \prod_{n=1}^\infty \left( 1 + \sum_{k=1}^\infty 3 \binom{k+2}{2}q^{kn} \right).$$

Portions of the argument are similar.  Again, we only need to look at parts of 1 or 2 part sizes.

Modulo 27, the contribution of $3\binom{k+2}{2}$ to the coefficient on $q^n$ is the number of partitions of $n$ satisfying the following conditions, multiplied by associated factors:

\begin{itemize}
\item exactly 2 part sizes both repeated a multiple of 3 times, multiplied by 9;
\item exactly 1 part size repeated 1 or 5 mod 9 times, multiplied by 9;
\item exactly 1 part size repeated 2 or 4 mod 9 times, multiplied by 18, and
\item exactly 1 part size repeated a multiple of 3 times, multiplied by 3.
\end{itemize}

But for $n \equiv 2 \pmod{9}$, the first and fourth are empty sets, and the other two divisors pair up as in the previous argument, and thus the congruence holds.

For $j=6$, we have 

\begin{multline*}C_{9,6}(q) \equiv_{27} \prod \left( 1+\frac{9q^n-9q^{2n}+30q^{3n}}{(1-q^n)^6} \right) \\ \equiv_{27} \prod \left( 1 + 9 \sum_{k=1}^\infty \binom{k+3}{4}q^{kn} + 3 \sum_{k=1}^\infty \binom{k+2}{5}q^{kn} \right).
\end{multline*}

Any partition counted with the $9\binom{k+3}{4}$ term contributing can only have 1 part size.  The binomial $\binom{k+3}{4}$ is 1 mod 3 if $k$ is 2 or 8 mod 9, and 2 mod 3 if $k$ is 1 or 7 mod 9; these divisors pair up as before.  For those with $3\binom{k+2}{5}$ contributing, see the previous discussion once again.

Finally, for $j \in \{1,4,7\}$, calculation of early cases suffices to eliminate many small $\ell$; for the case $j=1$ specifically, it suffices to note that $c_{9,1}(2) = 18 \not\equiv 0 \pmod{27}$, and hence no arithmetic progression starting with 2 will be 0 mod 27.  This concludes the proof of the theorem. \hfill $\Box$ 

\phantom{.}

\noindent \textbf{Remark:} Nicolas Smoot, in a preprint at \cite{Smoot}, has implemented an algorithm of Cristian-Silviu Radu for the verification of certain $\eta$-quotient identities, which effectively can take an arithmetic progression such as $3n+2$, a proposed modulus, and a generating function such as those for $j=8$ and $j=9$, and return a confirmation of whether the proposed congruence exists.  Upon presentation of this theorem at the Berndt Conference, he was able to verify the $j=8$ congruence with a few minutes of calculation.  Although it is not applicable to the other congruences listed here, the algorithm certainly bodes likely to be useful in the field more generally.

\subsection{Proof of Theorem \ref{2pp}}

Among the families of $C_{k,j}$ that generalize overpartitions, two that were studied in \cite{Keith1} were $C_{k,1}$ and $C_{k,k-1}$, both of which turn out to be $\eta$-quotients.  Another family that we might consider is $C_{2j,j}$.  For ease of calculation, we restrict ourselves to the case where $j$ is a prime $p$.

If we start with the basic generating function, we get $$C_{2p,p} (q) = \frac{\prod_{n=1}^\infty \left(1 + \binom{p}{1} q^n + \binom{p+1}{2} q^{2n} + \dots + \binom{2p-1}{p} q^{pn} \right)}{{(q;q)_\infty}^p}.$$

It is immediate that, when $p$ is prime, we get

$$C_{2p,p}(q) \equiv_p \prod_{n=1}^\infty \frac{1+q^{p n}}{1-q^{p n}}$$

\noindent which is just the overpartition function, magnified by $p$.  This easily gives us the relations 
$$ p \vert n \implies c_{2p,p}(n) \equiv \overline{p}(n/p) \pmod{p} \, , \quad \, \quad \text{otherwise} \, \, c_{2p,p}(n) \equiv 0 \pmod{p}.$$

But we can say a bit more mod $p^2$, over which $C_{2p,p}$ reduces to

$$C_{2p,p} \equiv_{p^2} \prod_{n=1}^\infty \frac{1+pq^n+2^{-1}pq^{2t}+\dots+(p-1)^{-1}pq^{(p-1)n}+q^{pn}}{1-pq^n-2^{-1}pq^{2t}-\dots-(p-1)^{-1}pq^{(p-1)n}-q^{pn}}$$

\noindent where inverses are taken mod $p^2$.

Each factor expands as $$1+ \sum_{t=1}^\infty 2 q^{p t n} + \sum_{t=0}^\infty (2t+2)p(1^{-1}q^{(tp+1)n} + \dots + (p-1)^{-1} q^{(tp+(p-1))n}).$$

This means that every part size contributes a factor of 2 to a weighting, and $2^{\# \{ \text{part sizes of } \, \lambda \} }$ is exactly the number of overpartitions of shape $\lambda$.  Hence partitions in which all part sizes are repeated a number of times divisible by $p$, which are in bijection with the partitions of $n/p$ by dividing the multiplicities, contribute a total of $\overline{p}(n/p)$ to the coefficient $c_{2p,p}(n)$ mod $p^2$.

If a part is repeated $tp+y$ times, $y \not\equiv 0 \pmod{p}$, such parts give an additional factor of $p(t+1)y^{-1}$ to the weight of the partition in the count.  Hence only partitions which have at most one such part can contribute a nonzero amount mod $p^2$.

Let $S$ be the set of $\lambda$ in which at most one part size is repeated a number of times not divisible by $p$, and let $t(\lambda)$ and $y(\lambda)$ be the associated values for such partitions.  Say $r(\lambda)$ is the number of part sizes of $\lambda$.  Taking the convention that the number of partitions, or overpartitions, of a non-integer amount is zero, we can refine the earlier result thus:

$$c_{2p,p}(n) \equiv_{p^2} \overline{p} \left( \frac{n}{p} \right) + p \sum_{\lambda \in S} 2^{r(\lambda)} (t(\lambda)+1) y(\lambda)^{-1} .$$

The theorem is proved. \hfill $\Box$

\phantom{.}

\noindent \textbf{Remark:} By isolating the part size repeated a nonzero number of times modulo $p$, one could probably further refine this theorem with some additional information to get a more elegant formula.

\subsection{Proof of Theorems \ref{K1Box} and \ref{SizesInBox}}

Each $r$ term on the right hand side of Theorem \ref{K1Box} is counting the partitions in which exactly $r$ part sizes appear; the multinomial term $(x_1 + \dots + x_k)^r$ distributes the colors among the part sizes in the usual way.  The meat of the theorem, then, is in Theorem \ref{SizesInBox}.

To prove the latter, what we need is the enumeration of overpartitions in the $M \times N$ box given by Jehan Dousse and Byungchan Kim in \cite{DoKim}: letting $s(\overline{\lambda})$ be the number of overlined parts of an overpartition $\overline{\lambda}$, they give

$$\sum_{\overline{\lambda} \in M \times N} t^{s(\overline{\lambda})} q^{\vert \overline{\lambda} \vert} = \sum_{k=0}^{min(M,N)} t^k q^{\binom{k+1}{2}} \left[ { M \atop k} \right]_q \left[ {{M+N-k} \atop N-k} \right]_q.$$

From our viewpoint, this is just setting $x_1 = t$, $x_2 = 1$ in Theorem \ref{K1Box}, where this identity becomes

$$\sum_{r=0}^{min(M,N)} (t+1)^r f_{(M \times N ; r)}(q) = \sum_{k=0}^{min(M,N)} t^k q^{\binom{k+1}{2}} \left[ { M \atop k} \right]_q \left[ {{M+N-k} \atop N-k} \right]_q.$$

Making the substitution $t \rightarrow t-1$, we have 

$$\sum_{r=0}^{min(M,N)} t^r f_{(M \times N ; r)}(q) = \sum_{k=0}^{min(M,N)} (t-1)^k q^{\binom{k+1}{2}} \left[ { M \atop k} \right]_q \left[ {{M+N-k} \atop N-k} \right]_q.$$

Extract coefficients on $t^r$ and we have the claim of the theorem. Apply this formula to Theorem \ref{K1Box} to get the explicit generating function.  \hfill $\Box$

\section{Future Directions}

The congruences described in Theorem \ref{mod9} are certainly just one example of a large family.  The study of congruences in $\eta$-quotients, as suggested by the remark concerning the work of Radu, has reached a high peak concerning the verification of any particular desired congruence, and much is known about congruences for overpartitions specifically, including larger infinite families.  Since $(k,j)$-colored partitions were originally proposed as a generalization of overpartitions, further developing the connections between these and the better-known overpartitions could be a useful line of investigation.

Still, most of the functions $C_{k,j}$ are not $\eta$-quotients, and so different tools will be necessary to obtain information about their congruences.  Where $\eta$-quotients have the power of well-developed modular form machinery that can be brought to bear, it seems that more combinatorial thinking may be fruitful in producing more information about $(k,j)$-colored partitions.  This could be a productive field of advancement.

An interesting conjecture arises when considering the coefficients in Theorem \ref{K1Box}.  Dousse and Kim conjecture that the coefficient of $q^n$ in $$\sum_{k=0}^{min(M,N)} t^k q^{\binom{k+1}{2}} \left[ { M \atop k} \right]_q \left[ {{M+N-k} \atop N-k} \right]_q$$ is unimodal in $t$.

This appears to generalize to more colors: if $x_k$ is set to 1, the array of the remaining coefficients appears to be unimodal in any linear cut.

For example, consider the $11 \times 7$ box with three colors, and look at the polynomial in $x_1$, $x_2$, $x_3$ which is the coefficient of $q^n$, $n=14$.  Set $x_3 = 1$.  The coefficients of $x_1^r x_2^c$, indexing row and column in order, are $$\begin{matrix} 101 & 291 & 300 & 129 & 19 \\ 291 & 600 & 387 & 76 & 0 \\ 300 & 387 & 114 & 0 & 0 \\ 129 & 76 & 0 & 0 & 0 \\ 19 & 0 & 0 & 0 & 0 \end{matrix}.$$

The antidiagonals are of course Pascal's triangle multiplied by the entries of the first column, as is immediate from Theorem \ref{K1Box}.  The rows (and columns, by symmetry) and diagonals also appear to be unimodal.  Even slant cuts through the array -- here, the sequence $(300, 387, 19)$ is the $(2,1)$ cut starting from position $(3,1)$ -- appear to have this property.  Perhaps Dousse and Kim's conjecture could be fruitfully set in this larger context.

\end{document}